\documentclass[12pt]{amsart}
\usepackage{amsmath}
\usepackage{enumerate}
\usepackage{amssymb}
\usepackage{amsbsy}
\usepackage{amsfonts}
\usepackage[arrow]{xy}

\allowdisplaybreaks
\theoremstyle{plain}

\newtheorem{thm}{Theorem}[section]
\newtheorem{prop}[thm]{Proposition}
\newtheorem{lem}[thm]{Lemma}
\newtheorem{cor}[thm]{Corollary}
\theoremstyle{definition}

\newtheorem{rmk}[thm]{Remark}
\newtheorem*{assum}{Assumption}
\numberwithin{equation}{section}
\newcommand{\sm}{\left(\begin{smallmatrix}}
\newcommand{\esm}{\end{smallmatrix}\right)}

\newfont{\FieldFont}{msbm10 scaled\magstep1}

\newcommand{\pf}{\noindent\bf Proof }

\def\re{\textrm{Re}}
\def\im{\textrm{Im}}

\usepackage{color}

\definecolor{blue}{rgb}{0,0,1}
\definecolor{red}{rgb}{1,0,0}
\definecolor{green}{rgb}{0,.6,.2}
\definecolor{purple}{rgb}{1,0,1}

\long\def\red#1\endred{{\color{red}#1}}
\long\def\blue#1\endblue{{\color{blue}#1}}
\long\def\purple#1\endpurple{{\color{purple}#1}}
\long\def\green#1\endgreen{{\color{green}#1}}

\begin{document}

\title[  Simultaneous nonvanishing of products of $L$-functions ]{ Simultaneous nonvanishing of products of $L$-functions associated to elliptic cusp forms}

 \author{YoungJu  Choie}

 \address{Department of Mathematics  \\
 Pohang University of Science and Technology\\
 Pohang, 790--784, Korea}
  \email{yjc@postech.ac.kr}
  
  \author{Winfried Kohnen}
\address{Department of Mathematics, University of Heidelberg, Germany  }
\email{winfried@mathi.uni-heidelberg.de}

\author{Yichao Zhang}

 \address{Department of Mathematics and Institute of Advanced Studies of Mathematics\\
 Harbin Institute of Technology\\ Harbin, 150001, P.R.China
  }
  \email{yichao.zhang@hit.edu.cn}

\date{\today}

\thanks{Key words: elliptic cusp forms, products of  $L$-functions, kernel function, double Eisenstein series}
 \thanks{2010
 Mathematics Subject Classification: primary 11F11, 11F67, 11F30; secondary  11F99}
 % \thanks{}

\begin{abstract} A generalized Riemann hypothesis  states that all zeros of the completed Hecke $L$-function  $L^*(f,s)$ of a normalized Hecke eigenform $f$ on the full modular group  should lie on the vertical line $Re(s)=\frac{k}{2}.$ 
It was shown in \cite{K} that there exists a Hecke eigenform $f$ of weight $k$ such that $L^*(f,s) \neq 0$ for sufficiently large $k$ and any point on the line segments $Im(s)=t_0, \frac{k-1}{2} < Re(s) < \frac{k}{2}-\epsilon, \frac{k }{2}+\epsilon < Re(s) < \frac{k+1}{2},$ for any given  real number $t_0$ and a positive real number $\epsilon.$
 This paper concerns   the non-vanishing of  
 the product $L^*(f,s)L^*(f,w)$ $(s,w\in \mathbb{C})$  on average. 

\end{abstract}
 \maketitle

%%%%%%%%%%%%%%%%%%%%%%%%%%%%%%%%%%%%%%%%%%%%%
%%%%%%%%%%%%%%%%%%%%%%%%%%%%%%%%%%%%%%%%%%%%%%%

\section{{\bf Introduction}}
Let $L^*(f,s)\, (s\in \mathbb{C})\, $ be the complete Hecke $L$-function of a non-zero cuspidal Hecke eigenform $f$ of integral weight $k$ on $SL_2(\mathbb{Z}).$  Although  the generalized Riemann hypothesis, which  states that all zeros should lie on the vertical line $Re(s)=\frac{k}{2},$ seems too far away  to prove at  this stage,  it is well known that zeros of $L^*(f,s)$ can occur only inside the critical strip 
 $\frac{k-1}{2} < Re(s) < \frac{k}{2}.$ However, Kohnen ( \cite{K}) showed that  for $k$ sufficiently large there exists a Hecke eigenform $f$ of weight $k$ such that $L^*(f,s) \neq 0$ at any point on the line segments $Im(s)=t_0, \frac{k-1}{2} < Re(s) < \frac{k}{2}-\epsilon, \frac{k }{2}+\epsilon < Re(s) < \frac{k+1}{2},$ for any given  real number $t_0$ and a positive real number $\epsilon.$
  This result and its method inspired various works on non-vanishing of L-values for different kinds of modular forms  (see  \cite {DK, HR, MP}).
\medskip

This paper concerns   the non-vanishing of  the product
 $L^*(f,s)L^*(f,w)$ $(s,w\in \mathbb{C})$ on average. 
 We shall prove that, given positive real numbers $T$ and $\delta$ and for all $k$ large enough  
  the sum of the products $L^*(f,s)L^*(f,w)$   over the   basis of Hecke eigenforms of weight $k$ does not vanish  
  on the region $Im(s), Im(w)\in [-T,\, T], \frac{k-1}{2}< Re(s), Re(w) <\frac{k+1}{2}, |Re(s)-\frac{k}{2}| >\delta, |Re(w)-\frac{k}{2}| >\delta.$ 
For the proof we compute the Fourier coefficients of  the double Eisenstein series, which is of independent interest.  Since it is   
 dual with respect to the  Petersson  scalar product of the values $L^*(f,s)L^*(f,w) $  by \cite{CD2}, 
 we derive the result by  estimating  the first term of Fourier coefficients.   This seems the first  non-vanishing  result for the product   $L^*(f,s)L^*(f,w) $ inside  the critical region. 
\bigskip

%%%%%%%%%%%%%%%%%%%%%%%%%%%%%%%%%%%%%%%%%%%
\section{{\bf{Notation}}}
\begin{itemize}
\item  For   complex numbers $z $ and $s$ with $z\neq 0$, fix the branch of  $z^s=e^{s\log{z}}$ as  $\log{z}=\log|z|+i \arg(z)$ and $-\pi < \arg(z) \leq \pi.$ 
\item $\mathbb{H} :$  the complex upper half plane. 
\item  $k :$ an even positive integer. 
\item $\Gamma :$   the full modular group.
\item  $\Gamma_{\infty}: $ the subgroup  generated by $\sm 1 & 1 \\ 0 & 1\esm.$
\item    $ \gamma=\sm a_{\gamma}& b_{\gamma} \\ c_{\gamma} & d_{\gamma} \esm$: the typical notation we employ for entries of a matrix $\gamma.$  

%$a_\gamma, c_\gamma$: the left-upper entry and left-lower entry of a %$2\times 2$ matrix $\gamma$ respectively.  

\item   $\mathcal{M}_n:=\left\{  \gamma \in M_{2\times 2}(\mathbb{Z})  : \, \det(\gamma)=n\right\}.$
 \item $\zeta(s)=\sum_{n\geq 1}\frac{1}{n^s}$:   the Riemann zeta function.
\item    $S_k$   the space of cusp forms of weight $k$ on $\Gamma $ with the Petersson scalar product $<\, \, , >.$  
\item  $\mathcal{B}_k$ the basis of normalized Hecke eigenforms, i.e,  eigenforms whose first Fourier coefficient equals $1$, of $S_k$.

\end{itemize}

\bigskip
%%%%%%%%%%%%%%%%%%%%%%%%%%%%%%%%%%%%%%%%
\section{{\bf Statement of Result}}
For $f(\tau)=\sum_{n\geq 1} a_f(n)e^{2\pi i n \tau}$ $(\tau \in \mathbb{H}),$ let $L(f,s):=\sum_{n\geq 1} \frac{a_f(n)}{n^s}$ be the Hecke $L$-function associated to $f.$  It is well-known that the  completed $L$-function 
$$L^*(f,s):=(2\pi)^s\Gamma(s)L(f,s)$$ has an analytic continuation and the functional equation
$$L^*(f, k-s)=(-1)^{\frac{k}{2}}L^*(f,s).$$

\begin{thm}[Main Theorem]\label{main-thm}
For any  fixed positive real numbers $T,\delta$, let the region $\mathcal{R}_{T,\delta}$ of points $(s,w)\in\mathbb{C}^2$ such that
\begin{itemize}
\item $\im(s),\im(w)\in[-T,T]$, 
\item $\frac{k-1}{2}<\re(s),\re(w)<\frac{k+1}{2}$, 
\item $|\re(s)-\frac{k}{2})|\geq \delta$ and $|\re(w)-\frac{k}{2})|\geq \delta$.
\end{itemize}
Then there exists a   constant $C(T,\delta) >0 $  depending only on $T$ and $\delta$ such that for 
 $k>C(T,\delta)$, the following function 
\[\sum_{f\in \mathcal{B}_k}\frac{L^*(f,s )L^*(f,w )}{\langle f,f\rangle} \]
does  not vanish 
at  any pair $(s,w)\in\mathcal{R}_{T,\delta}.$
\end{thm}

\begin{cor}
Let $T$ and $\delta$ be positive real numbers.  Then for $k > C(T,\delta) $ and any pair of complex numbers $(s,w)\in \mathcal{R}_{T,\delta},$ there exists a Hecke eigenform $f\in S_k$ such that  $L^*(f,s )L^*(f,w ) \neq 0.$
\end{cor}

\medskip

On the region $ \mathcal{R}_{T,\delta}$ in Theorem \ref{main-thm},  $s,w$ are away from the central lines $\re(s),\re(w)=\frac{k}{2}.$   The points with $s= \frac{k}{2}$ or $w= \frac{k}{2}$ have to be removed, since the L-values $L^*(f,\frac{k}{2})$ are necessarily $0 $ for odd $\frac{k}{2} .$
However we may try to  enlarge the non-vanishing region by adding  points with $\re(s)$ and $\re(w)$  equal to $\frac{k}{2}. $    It turns out that this is closely related to   a certain  property of the Riemann zeta function $\zeta(s)$. Since we do not know whether such property holds in general, we add it as an assumption. 

\begin{assum}\label{ass} Let $z_0\in\mathbb{C}$ with $\re(z_0)>0$ and $\im(z_0)\neq 0$, and $s_0=it_0$ with $t_0\neq 0$. Then there exist $\epsilon>0$, a neighborhood $W\subset\mathbb{C}^2$ of $(s_0,z_0)$ and a positive integer $N$, such that if $k\geq N$ and $(s,z)\in W$ then
\[\left|\left(\frac{4\pi}{k}\right)^s\frac{\zeta(1+z+s)}{\zeta(1+z-s)}\right|>1+\epsilon \textrm{ or } <1-\epsilon.\]
\end{assum}
\medskip

Note that when $W$ is sufficiently small, the arguments in the zeta values stay in the right half plane $\re(w)>1$, which makes the zeta values non-vanishing.   
%The Assumption is simply one that can make the proof of the following theorem to
%work.  
Roughly speaking, the above assumption says that the magnitude of the Riemann zeta function $|\zeta(s)|$ is not (locally) symmetric about a non-real horizontal line on the right-half plane $\re(s)>1$, so in particular it is likely not related to the Riemann Hypothesis. 

\begin{thm}\label{without} For any  fixed positive real numbers $T,\delta$, let the region  $\widetilde{\mathcal{R}_{T,\delta}}$  of points $(s,w)\in\mathbb{C}^2$ such that
\begin{itemize}
\item $\delta\leq |\im(s)|,|\im(w)|\leq T$, 
\item $\frac{k-1}{2} <\re(s),\re(w)<\frac{k+1}{2}$, 
\item   $|\re(s)-\frac{k}{2}|+|\re(w)-\frac{k}{2}|\geq \delta$,
\end{itemize}
Suppose the preceding assumption holds. Then there exists a   constant $C(T,\delta) >0 $  depending only on $T$ and $\delta$ such that for $k>C(T,\delta),$
\[\sum_{f\in \mathcal{B}_k}\frac{L^*(f,s)L^*(f,w)}{\langle f,f\rangle} \]
does not vanish at any pair $(s,w) \in \widetilde{\mathcal{R}_{T,\delta}} .$
\end{thm}

\begin{rmk} The condition $ |\re(s)-\frac{k}{2}|+|\re(w)-\frac{k}{2}|\geq \delta$
in Theorem \ref{without} looks striking since   $Re(s)=\frac{k}{2}$ may be chosen and the generalized Riemann hypothesis states that the zeros should lie on this line.  However  the implication of Theorem \ref{without} shows that 
for any fixed nonzero real number $t_0$, there exists $C(t_0)>0$ such that whenever $k>C(t_0)$ there exists $f\in\mathcal{B}_k$ with $L^*(f,\frac{k}{2}+it_0)\neq 0$.
 \end{rmk}
 
\medskip

\noindent {\bf Acknowledgment}: The first author was supported by NRF2018R1A4A1023590 and \\
NRF2017R1A2B2001807. The third author was partially supported by NSFC11871175.   After finishing this paper, we learned that the method given in \cite{BF} may be modified to cover our situation.   

The authors would like to thank the referee  for  helpful comments and suggestions which greatly improved the exposition of this paper

\bigskip

 %%%%%%%%%%%%%%%%%%%%%%%%%%%%%%%%%%%%%%%%%%%%%%%%%%%%%%%%%%%%%%%%%%%%%%
\section{\bf{Proof}}
 
%%%%%%%%%%%%%%%%%%%%%%%%%%%%%%%%%%%%%%%%%%%%

\subsection{Double Eisenstein Series and its Fourier expansion} 
We recall the basics on double Eisenstein series (see \cite{CD2, CD1}), and following the lines in \cite{K}, we compute its Fourier expansion.
\medskip

For $s\in \mathbb{C}$ the double Eisenstein series is defined as
\[E_{s,k-s}(z,w)=\sum_{\gamma,\delta\in \Gamma_{\infty}\backslash\Gamma, c_{\gamma\delta^{-1}}>0}(c_{\gamma\delta^{-1}})^{w-1}\left(\frac{j(\gamma,z)}{j(\delta,z)}\right)^{-s}j(\delta,z)^{-k}.\]
In \cite{CD2} it is shown that $E_{s,k-s}(z,w)$ converges absolutely and uniformly on compact sets for $\mathcal{D},  $ where
\begin{equation} \label{reg-1}
\mathcal{D}:=\{ 2<\re(s)<k-2,\quad \re(w)<\min\left\{\re(s)-1,k-\re(s)-1\right\}\}.
\end{equation}
 
Define the completed double Eisenstein series
\[E_{s,k-s}^*(z,w)=\frac{e^{\pi is/2}\Gamma(s)\Gamma(k-s)\Gamma(k-w)\zeta(1-w+s)\zeta(1-w+k-s)}{2^{3-w}\pi^{k+1-w}\Gamma(k-1)}E_{s,k-s}(z,w).\] Then the following  are proven :
\begin{thm}\cite{CD2}
Let $k\geq 6$ be even. The series
$E_{s,k-s}^*(z, w)$ has an analytic continuation to all $s,w \in \mathbb{C}$ and
as a function of $z$ is always in $S_k.$
For any $f\in \mathcal{B}_k,$ we have
$$<E_{s,k-s}^*(z, w)\, , f>=L^*(f,s)L^*(f,w).$$
\end{thm}

 The following are consequences of the above results given in \cite{CD2}:
 
\begin{itemize}
\item\label{Eisen-functional}
The functional equations of $L^*(f,s)$ induces those of $E^*_{s,k-s}(z,w) :$
 $$E^*_{s,k-s}(z,w) = E^*_{w,k-w}(z,s),\quad 
 E^*_{s,k-s}(z,w)= (-1)^{\frac{k}{2}} E^*_{k-s, s}(z,w).$$ 
 \item
 \label{Prod}
\begin{eqnarray*}
E_{s,k-s}^*(z, w)&=&\sum_{f\in \mathcal{B}_k}\frac{L^*(f,s)L^*(f,w)}{\langle f,f\rangle}f(z)
 \nonumber\\
\label{main}
&=&\frac{e^{\pi is/2}\Gamma(s)\Gamma(k-s)\Gamma(k-w)}{2^{3-w}\pi^{k+1-w}\Gamma(k-1)}\sum_{n=1}^\infty n^{w-1}\sum_{\gamma\in\mathcal{M}_n}(\gamma z)^{-s}j(\gamma,z)^{-k} 
\end{eqnarray*}
so that the right-hand side has analytic continuation to all of $(s,w)\in\mathbb{C}^2.$
\end{itemize}

\begin{prop}({\bf Fourier expansion}) \label{Fourier}   Let $(s,w)\in \mathcal{D}$  as in \eqref{reg-1}.   
$E^*_{s,k-s}(z,w)$ has the Fourier expansion
\[E^*_{s,k-s}(z,w)=\frac{\Gamma(s)\Gamma(k-s)\Gamma(k-w)}{2^{2-w}\pi^{k+1-w}\Gamma(k-1)}\sum_{m=1}^\infty c_{s,w,k}(m)e^{2\pi imz},\] where
\begin{align*}
&c_{s,w,k}(m)\\
=& \frac{(2\pi)^s}{\Gamma(s)}m^{s-1}\sigma_{w-s}(m)\zeta(k-s-w+1)\\
&+(-1)^\frac{k}{2}\frac{(2\pi)^{k-s}}{\Gamma(k-s)}m^{k-s-1}\sigma_{w+s-k}(m)\zeta(s-w+1)\\
&+(-1)^{\frac{k}{2}}\frac{(2\pi)^km^{k-1}}{\Gamma(s)\Gamma(k-s)}\sum_{a,c>0,(a,c)=1}c^{s-k}a^{-s}\sum_{n\geq 1}n^{w-1}\sum_{r\mid m}r^{w-k}\\&\ \times\left(e^{\pi is/2}e^{2\pi i \frac{m}{r}\frac{na'}{c}}\sideset{_1}{_1}{\mathop{f}}\left(s,k;-\frac{2\pi imn}{rac}\right)+e^{-\pi is/2}e^{-2\pi i \frac{m}{r}\frac{na'}{c}}\sideset{_1}{_1}{\mathop{f}}\left(s,k;\frac{2\pi imn}{rac}\right)\right)
\end{align*} 
  is absolutely convergent. Here 
\[\sideset{_1}{_1}{\mathop{f}}\left(\alpha,\beta;z\right)=\frac{\Gamma(\alpha)\Gamma(\beta-\alpha)}{\Gamma(\beta)}\sideset{_1}{_1}{\mathop{F}}\left(\alpha,\beta;z\right) \] with Kummer's degenerate hypergeometric function ${}_1F_1(\alpha, \beta; z)$  and $a'=a^{-1} \pmod{c}.$ 
\end{prop}

%\blue We need that $(s,w)\in\mathcal{D}$ to show the absolute convergence of the series for %the third term of the Fourier coefficient. For the other terms, it does not matter, since they are all %meromorphic functions. Maybe we do not revise it but explain it in the letter instead?\endblue

{\pf}of Proposition \ref{Fourier}:
 To compute the Fourier expansion of 
\[ \sum_{n=1}^\infty n^{w-1}\sum_{\gamma\in\mathcal{M}_n}(\gamma z)^{-s}j(\gamma,z)^{-k}=\sum_{m=1}c_{s,w,k}(m)e^{2\pi i m z} \] 
we split it into four cases:\\

\noindent ({\bf{1}}) Consider first  the elements $\gamma\in\mathcal{M}_n$ with $c_\gamma=0$. The contribution of such terms to the $m$-th Fourier coefficient is given by
\begin{align*}
{\bf I}:=&\sum_{n\geq 1}n^{w-1}\sum_{ad=n}\sum_{b\in\mathbb{Z}}\int_{iC}^{iC+1}\left(\frac{az+b}{d}\right)^{-s}d^{-k}e^{-2\pi i mz}dz\\
=& 2 \sum_{n\geq 1}n^{w-1}\sum_{ad=n,a>0}\sum_{b\in\mathbb{Z}}\int_{iC}^{iC+1}\left(\frac{az+b}{d}\right)^{-s}d^{-k}e^{-2\pi i mz}dz\\
=& 2 \sum_{n\geq 1}n^{w-1}\sum_{ad=n,a>0}d^{s-k}\int_{iC}^{iC+1}\sum_{b\in\mathbb{Z}}\left(az+b\right)^{-s}e^{-2\pi i mz}dz,
\end{align*}
where $C$ is any fixed positive real number. Note that since we only work on $\mathcal{D}$ in this section, all interchanges of sums and/or integrals are justified by the absolute convergence.
Applying Lipschitz's formula
\begin{equation}\label{Lip}
\sum_{n\in\mathbb{Z}}(z+n)^{-s}=\frac{e^{-\pi is/2}(2\pi)^s}{\Gamma(s)}\sum_{n\geq 1}n^{s-1}e^{2\pi in\tau}, \quad \im(\tau)>0, \re(s)>1,
\end{equation}
to the sum over $b$, we have
\begin{eqnarray}\label{C-1}
{\bf I} &=& 2 \frac{e^{-\pi is/2}(2\pi)^s}{\Gamma(s)}\sum_{n\geq 1}n^{w-1}\sum_{ad=n,a>0}d^{s-k}\int_{iC}^{iC+1}\sum_{r\geq 1}r^{s-1}e^{2\pi i raz}e^{-2\pi i mz}dz\nonumber \\
&=& 2 \frac{e^{-\pi is/2}(2\pi)^s}{\Gamma(s)}\sum_{n\geq 1}n^{w-1}\sum_{ad=n,a\mid m}d^{s-k}(m/a)^{s-1}.
\nonumber \\
&=& 2 \frac{e^{-\pi is/2}(2\pi)^s}{\Gamma(s)}\sum_{a\mid m}\sum_{d\geq 1} (ad)^{w-1}d^{s-k}(m/a)^{s-1},  
\nonumber \\
&=& 2 \frac{e^{-\pi is/2}(2\pi)^s}{\Gamma(s)}  m^{s-1}\sigma_{w-s}(m)\zeta(k-s-w+1),   
\end{eqnarray}
where $\sigma_s(n)=\sum_{d\mid n}d^s$ is the divisor function.
\medskip

\noindent  ({\bf 2}) Similarly, we obtain the contribution {\bf II}  of the terms $\gamma\in\mathcal{M}_n$ with $a_\gamma=0$ in the $m$-th Fourier coefficient:
\begin{equation}\label{C-2}
{\bf II}=2 \frac{e^{\pi i(k-s)/2}(2\pi)^{k-s}}{\Gamma(k-s)}m^{k-s-1}\sigma_{w+s-k}(m)\zeta(s-w+1).
\end{equation}
  Note the symmetry $s\rightarrow k-s$ when switching from $\bf{I}$ to ${\bf{II}}.$  
\medskip

\noindent  ({\bf 3}) Next we consider the contribution {\bf III} of the terms $\gamma\in\mathcal{M}_n$ with $a_\gamma c_\gamma > 0$ in the $m$-th Fourier coefficient. 
 The set of integral matrices with determinant $n$ can be listed as follows:
\[\left\{
\begin{pmatrix}
ar & nb_0/r+(t+r\ell)a\\
cr & nd_0/r+(t+r\ell)c
\end{pmatrix}\colon r\mid n, \textrm{gcd}(a,c)=1, t\in \mathbb{Z}/r\mathbb{Z},\ell\in\mathbb{Z}
\right\}\]
where for each pair $(a,c)$, $b_0,d_0$ are fixed so that $ad_0-b_0c=1$. With this, we have

\begin{eqnarray}\label{ee}
{\bf III} &=&  \sum_{n\geq 1}n^{w-1}\sum_{ac>0,(a,c)=1}\sum_{r\mid n,t\in\mathbb{Z}/r\mathbb{Z}}\sum_{\ell\in\mathbb{Z}} 
 \nonumber\\
&& \times\int_{iC}^{iC+1}\left(\frac{ra(z+\ell)+nb_0/r+ta}{rc(z+\ell)+nd_0/r+tc}\right)^{-s}(rc(z+\ell)+nd_0/r+tc)^{-k}e^{-2\pi i m(z+\ell)}dz \nonumber \\
&=&\sum_{n\geq 1}n^{w-1}\sum_{ac>0,(a,c)=1}\sum_{r\mid n,t\in\mathbb{Z}/r\mathbb{Z}}\nonumber \\
&&\times\int_{iC-\infty}^{iC+\infty}\left(\frac{raz+nb_0/r+ta}{rcz+nd_0/r+tc}\right)^{-s}(rcz+nd_0/r+tc)^{-k}e^{- 2\pi i m(z+\ell)}dz 
\nonumber \\
 && \mbox{ (since $e^{2\pi i m \ell}=1$ and by the change of variable $z\rightarrow z-(nd_0+trc)/(r^2c)$) }\nonumber\\
&=& \sum_{n\geq 1}n^{w-1}\sum_{ac>0,(a,c)=1}\sum_{r\mid n,t\in\mathbb{Z}/r\mathbb{Z}}e^{2\pi i m\frac{nd_0+trc}{r^2c}}\nonumber \\
&&\times\int_{iC-\infty}^{iC+\infty}\left(\frac{a}{c}-\frac{n}{c^2r^2z}\right)^{-s}(rcz)^{-k}e^{-2\pi i mz}dz.
\nonumber \\
&& \mbox{(since the summation over $t$ vanishes unless $r\mid m$ )}\nonumber\\
&=& \sum_{n\geq 1}n^{w-1}\sum_{ac>0,(a,c)=1}\sum_{r\mid (n,m)}re^{2\pi i \frac{m}{r}\frac{n}{r}\frac{a'}{c}}\nonumber \\
&& \times\int_{iC-\infty}^{iC+\infty}\left(\frac{a}{c}-\frac{n}{c^2r^2z}\right)^{-s}(rcz)^{-k}e^{-2\pi i mz}dz 
 \label{Arg}\\
&& \mbox{(with   $0<a'\leq c$ and  $a'a \equiv 1\pmod{c}$ )}\nonumber\\
&=& \sum_{n\geq 1}n^{w-1}\sum_{ac>0,(a,c)=1}\sum_{r\mid (n,m)}re^{2\pi i \frac{m}{r}\frac{n}{r}\frac{a'}{c}}\nonumber\\
&&\times\int_{iC-\infty}^{iC+\infty}z^s\left(\frac{a}{c}z-\frac{n}{c^2r^2}\right)^{-s}(rcz)^{-k}e^{-2\pi i mz}dz  \nonumber\\
&=&\sum_{n\geq 1}n^{w-1}\sum_{ac>0,(a,c)=1}\sum_{r\mid (n,m)}r^{1-k}c^{-k}e^{2\pi i \frac{m}{r}\frac{n}{r}\frac{a'}{c}}\nonumber \\
&&\times\int_{iC-\infty}^{iC+\infty}\left(\frac{a}{c}z-\frac{n}{c^2r^2}\right)^{-s}z^{s-k}e^{-2\pi i mz}dz.\nonumber\\
&& \mbox{(the change of variable $z\mapsto \frac{c}{a}iz$)}\nonumber\\
&=& \sum_{n\geq 1}n^{w-1}\sum_{ac>0,(a,c)=1}\sum_{r\mid (n,m)}(-1)^{k/2}(c/a)^{s-k+1}r^{1-k}c^{-k}e^{2\pi i \frac{m}{r}\frac{n}{r}\frac{a'}{c}}\nonumber\\
&& \times \frac{1}{i}\int_{C-i\infty}^{C+i\infty}\left(z+\frac{in}{c^2r^2}\right)^{-s}z^{s-k}e^{2\pi mcz/a}dz \nonumber\\
 && \mbox{ $\bigl($ using the following integral representation }  \nonumber\\
  &&  
  \frac{1}{2\pi i}\int_{C-i\infty}^{C+i\infty}
(z+\alpha)^{-\mu}(z+\beta)^{-\nu}e^{pz}dz= 
 \frac{p^{\mu + \nu-1}e^{-\beta p} }{\Gamma(\mu+\nu)}
 \sideset{_1}{_1}{\mathop{F}}(\mu,\mu+\nu ; (\beta- \alpha)p),
  \nonumber\\
&& \mbox{  for all $p,\mu,\nu \in \mathbb{C} $ with  $\re(\mu),\re(\nu)>0 \bigr)$}\nonumber\\ 
 &=&  2(-1)^{\frac{k}{2}}\frac{(2\pi)^km^{k-1}}{\Gamma(k)}\sum_{a,c>0,(a,c)=1}c^{s-k}a^{-s}\sum_{n\geq 1}n^{w-1}\nonumber \\
&&\times\sum_{r\mid m}r^{w-k}e^{2\pi i \frac{m}{r}\frac{na'}{c}}\sideset{_1}{_1}{\mathop{F}}\left(s,k;-\frac{2\pi imn}{rac}\right). \label{C3+}
\end{eqnarray}

\noindent {\bf(4)}  Finally the computation on terms with $ac<0$ can be done similarly. One has to pay attention to the first equality in \eqref{Arg} and compute using
\[\left(\frac{a}{c}-\frac{n}{c^2r^2z}\right)^{-s}=(-1)^{-s}\left(-\frac{a}{c}+\frac{n}{c^2r^2z}\right)^{-s}=e^{-\pi is}z^{-s}\left(-\frac{a}{c}z+\frac{n}{c^2r^2}\right)^{-s}.\] Then by replacing $(a,c)$ with $(-a,c)$, the above computation in the case of $ac>0$ shows that
\begin{align}\label{C3-}
 {\bf IV} 
=&2(-1)^{\frac{k}{2}}\frac{(2\pi)^km^{k-1}}{e^{\pi is}\Gamma(k)}\sum_{a,c>0,(a,c)=1}c^{s-k}a^{-s}\sum_{n\geq 1}n^{w-1}\\&\times\sum_{r\mid m}r^{w-k}e^{-2\pi i \frac{m}{r}\frac{na'}{c}}\sideset{_1}{_1}{\mathop{F}}\left(s,k;\frac{2\pi imn}{rac}\right).\nonumber
\end{align}

Combining the formulas \eqref{C-1}, \eqref{C-2}, \eqref{C3+} and \eqref{C3-} we conclude the result.
{\qed}
\medskip

\begin{cor}\label{D}  The formula
\begin{align*} 
&\frac{2^{2-w}\pi^{k+1-w}\Gamma(k-1)}{\Gamma(s)\Gamma(k-s)\Gamma(k-w)}\sum_{f\in \mathcal{H}_k}\frac{L^*(f,s)L^*(f,w)}{\langle f,f\rangle}\\
=&\frac{(2\pi)^s}{\Gamma(s)}\zeta(k-s-w+1)+(-1)^\frac{k}{2}\frac{(2\pi)^{k-s}}{\Gamma(k-s)}\zeta(s-w+1)\\
&+(-1)^{\frac{k}{2}}\frac{(2\pi)^k}{\Gamma(s)\Gamma(k-s)}\sum_{a,c>0,(a,c)=1}c^{s-k}a^{-s}\sum_{n\geq 1}n^{w-1}\\&\ \times\left(e^{\pi is/2}e^{2\pi i na'/c}\sideset{_1}{_1}{\mathop{f}}\left(s,k;-\frac{2\pi in}{ac}\right)+e^{-\pi is/2}e^{-2\pi ina'/c}\sideset{_1}{_1}{\mathop{f}}\left(s,k;\frac{2\pi in}{ac}\right)\right)
\end{align*}
holds on $\mathcal{D}$ and the double sum in the last term of the right-hand side is absolutely convergent on $\mathcal{D}$.
\end{cor}
{\pf}of Corollary \ref{D}:
It follows easily from Proposition \ref{Fourier} and the formula \eqref{Prod}.
{\qed}

\begin{rmk} 
By (5-2) of \cite{CD2}, one may apply Kohnen's Fourier expansion for Cohen's kernel function (\cite{K}, Lemma 2) and obtain the same formula as in Corollary \ref{D}. However, 
the absolute convergence on $\mathcal{D}$ of the double sum in Corollary \ref{D} is no longer clear. In our approach, the absolute convergence of $E^*_{s,k-s}(z,w)$ implies that of the double sum together with the integral for the hypergeometric function $\sideset{_1}{_1}{\mathop{F}}$, hence is stronger than  that of the double sum in Corollary \ref{D}. 
\end{rmk}

\bigskip

\subsection{{\bf Analytic continuation}}The left-hand side of the identity in Corollary \ref{D} is meromorphic on $\mathbb{C}^2$, while the right-hand side is only valid on $\mathcal D$. For later purpose, we shall analytically continue the right-hand side, and to this end, the following domains  will be involved:
\begin{align*}
\mathcal{D}_1:=\{(s,w)\in \mathbb{C}^2:   &\quad 2<\re(s)<k-2,\quad \re(w)<0\},\\
\mathcal{F}:=\{(s,w)\in \mathbb{C}^2:  &\quad 3/2<\re(s),\re(w)<k-2\}.
\end{align*}

\begin{prop}\label{F} We have the following identity on $\mathcal{F}$
\begin{align*} 
&\frac{2^{2-w}\pi^{k+1-w}\Gamma(k-1)}{\Gamma(s)\Gamma(k-s)\Gamma(k-w)}\sum_{f\in \mathcal{H}_k}\frac{L^*(f,s)L^*(f,w)}{\langle f,f\rangle}\\
=&\frac{(2\pi)^s}{\Gamma(s)}\zeta(k-s-w+1)+(-1)^\frac{k}{2}\frac{(2\pi)^{k-s}}{\Gamma(k-s)}\zeta(s-w+1)\\
&+2(-1)^{\frac{k}{2}}\frac{(2\pi)^{k-w}\Gamma(w)\Gamma(s-w)}{\Gamma(s)\Gamma(k-w)}\cos(\pi(s-w)/2)\zeta(s-w)\\
&+2(-1)^{\frac{k}{2}}\frac{(2\pi)^{k-w}\Gamma(w)\Gamma(k-s-w)}{\Gamma(k-s)\Gamma(k-w)}\cos(\pi(s+w)/2)\\
&+ R(s,w),
\end{align*}
where $R(s,w)$ is holomorphic on $\mathcal{F}$ and bounded by \eqref{Remainder}.
\end{prop}
{\pf}of Proposition \ref{F}: 
 We only have to deal with the last term of the right-hand side of Corollary \ref{D}, namely
\begin{align*}A(s,w):&=(-1)^{\frac{k}{2}}\frac{(2\pi)^k}{\Gamma(s)\Gamma(k-s)}\sum_{a,c>0,(a,c)=1}c^{s-k}a^{-s}\sum_{n\geq 1}n^{w-1}\\&\ \times\left(e^{\pi is/2}e^{2\pi i na'/c}\sideset{_1}{_1}{\mathop{f}}\left(s,k;-\frac{2\pi in}{ac}\right)+e^{-\pi is/2}e^{-2\pi ina'/c}\sideset{_1}{_1}{\mathop{f}}\left(s,k;\frac{2\pi in}{ac}\right)\right).
\end{align*}
{{\bf Step I}}\, \, Let    $G_{a,c}(s,w):= G_{a,c}'(s,w) +G_{a,c}"(s,w)$,
where 
$$G_{a,c}'(s,w):=\sum_{n\geq 1}n^{w-1}  e^{\pi is/2}e^{2\pi i na'/c}\sideset{_1}{_1}{\mathop{f}}\left(s,k;-\frac{2\pi in}{ac}\right),
$$
$$G_{a,c}''(s,w):=\sum_{n\geq 1}n^{w-1}  e^{-\pi is/2}e^{-2\pi ina'/c}\sideset{_1}{_1}{\mathop{f}}\left(s,k;\frac{2\pi in}{ac}\right) 
$$
  As a subseries of an absolutely convergent series, the function $G_{a,c}(s,w)$ is holomorphic on $\mathcal{D}$, so in particular $G_{a,c}(s,w)$ is holomorphic on the smaller region $\mathcal{D}_1.$ 
For $\re(\beta)>\re(\alpha)>0$, it is known that (\cite{K}) 
\[\sideset{_1}{_1}{\mathop{f}}\left(\alpha,\beta;z\right)=\int_0^1e^{zu}u^{\alpha-1}(1-u)^{\beta-\alpha-1}du.\]
On $\mathcal{D}_1$, we have
\begin{align*}
G_{a,c}'(s,w)=&\sum_{n=1}^\infty n^{w-1}e^{2\pi ina'/c}e^{\pi is/2} \sideset{_1}{_1}{\mathop{f}}\left(s,k; \frac{-2\pi in}{ac}\right)\\
=&\sum_{n=1}^\infty n^{w-1}e^{2\pi ina'/c}e^{\pi is/2}\int_0^1e^{\frac{-2\pi inu}{ac}}u^{s-1}(1-u)^{k-s-1}du\\
=& e^{\pi is/2}\int_0^1u^{s-1}(1-u)^{k-s-1}\sum_{n=1}^\infty n^{w-1}e^{2\pi ina'/c}e^{\frac{-2\pi inu}{ac}}du,
\end{align*}
where the interchange of summation and integration is justified because of absolute convergence on $\mathcal{D}_1$. We warn here that the above series with $\sideset{_1}{_1}{\mathop{f}}$ expanded is not necessarily absolutely convergence on $\mathcal{D}$, which is why we need the smaller $\mathcal D_1$. Now for $s\in\mathbb C$, $a>0$, define 
\[F(s,a)=\sum_{n=1}^\infty\frac{e^{2\pi ina}}{n^s},\quad \zeta(s,a)=\sum_{n=0}^\infty\frac{1}{(n+a)^s}.\]
Therefore, originally on $\mathcal{D}_1$,
\begin{align*}
G_{a,c}'(s,w)&=e^{\pi is/2}\int_0^1u^{s-1}(1-u)^{k-s-1}F\left(1-w,\frac{a'}{c}-\frac{u}{ac}\right)du.
\end{align*}
  We specify $0< a'\leq c$ so that $0<\frac{a'}{c}-\frac{u}{ac}<1$. So the identity between $F(s,a)$ and $\zeta(s,a)$ (see Formula 25.13.2 of \cite{O})  implies that 
\begin{align*}
G_{a,c}'(s,w)=&(2\pi)^{-w}\Gamma(w)\int_0^1u^{s-1}(1-u)^{k-s-1}\\
&\times\left( e^{\pi i(s+w)/2}\zeta\left(w,\frac{a'}{c}-\frac{u}{ac}\right)+e^{\pi i(s-w)/2}\zeta\left(w,1-\frac{a'}{c}+\frac{u}{ac}\right)\right)du.
\end{align*}
 
  Similarly,
\begin{align*}
G_{a,c}''(s,w)&=e^{-\pi is/2}\int_0^1u^{s-1}(1-u)^{k-s-1}F\left(1-w,-\frac{a'}{c}+\frac{u}{ac}\right)du 
\end{align*}
and 
\begin{align*}
G_{a,c}''(s,w)=&(2\pi)^{-w}\Gamma(w)\int_0^1u^{s-1}(1-u)^{k-s-1}\\
&\times\left( e^{-\pi i(s+w)/2}\zeta\left(w,\frac{a'}{c}-\frac{u}{ac}\right)+e^{\pi i(w-s)/2}\zeta\left(w,1-\frac{a'}{c}+\frac{u}{ac}\right)\right)du.
\end{align*}

So  we have,  on $\mathcal{D}_1$ 
\begin{align}\label{G}
&G_{a,c}(s,w)\\
\nonumber =&2(2\pi)^{-w}\Gamma(w)\int_0^1u^{s-1}(1-u)^{k-s-1}\\
\nonumber&\times\left( \cos(\pi(s+w)/2)\zeta\left(w,\frac{a'}{c}-\frac{u}{ac}\right)+\cos(\pi(s-w)/2)\zeta\left(w,1-\frac{a'}{c}+\frac{u}{ac}\right)\right)du.
\end{align}
Note that the right-hand side of \eqref{G} is meromorphic on $\mathcal{D}$, forcing that \eqref{G} holds on $\mathcal{D}$ as well since $G_{a,c}(s,w)$ is holomorphic on $\mathcal{D}$.

 Replacing the expression $G_{a,c}(s,w)$ in $A(s,w)$ with the right-hand side of \eqref{G}, the following equality
\begin{align}\label{A}
A(s,w)&=2(-1)^{\frac{k}{2}}\frac{(2\pi)^{k-w}\Gamma(w)}{\Gamma(s)\Gamma(k-s)}\sum_{(a,c)=1,a,c>0}c^{-k+s}a^{-s}\int_0^1u^{s-1}(1-u)^{k-s-1}\\
\nonumber&\times \left( \cos(\pi(s+w)/2)\zeta\left(w,\frac{a'}{c}-\frac{u}{ac}\right)+\cos(\pi(s-w)/2)\zeta\left(w,1-\frac{a'}{c}+\frac{u}{ac}\right)\right)du
\end{align} 
holds on $\mathcal{D}$.

\noindent  {\bf Step II} \, Next, we prove that the series on the right-hand side of \eqref{A} is absolutely convergent on $\mathcal{F}.$
Since $\mathcal{D}\cap\mathcal{F}\neq\emptyset$, we then see that the equality \eqref{A} holds on $\mathcal{F}$.

\noindent({\bf  1}) First consider the terms with $c>1$. Note that for $\re(w)>\frac{3}{2}$ and each pair $a,c$, we always have $1-\frac{a'}{c}\geq \frac{1}{c}$, so 
\begin{align*}
&\left|\zeta\left(w,1-\frac{a'}{c}+\frac{u}{ac}\right)\right|\leq\left|\left(1-\frac{a'}{c}+\frac{u}{ac}\right)^{-w}\right|+\zeta(\re(w))\\
\leq &\left(\frac{1}{c}+\frac{u}{ac}\right)^{-\re(w)}+\zeta(3/2)\leq c^{\re(w)}+\zeta(3/2)\ll c^{\re(w)}.
\end{align*}
Then the sum over such pairs $(a,c)$ is bounded absolutely by   a constant multiple of 
\begin{align*}
&e^{\pi(|\im(s)|+|\im(w)|)}\sum_{c=2}^\infty\sum_{m=1}^\infty{\sum_{\substack{a\text{ mod} c\\ \gcd{(a,c)}=1}}}c^{-k+\re(s)+\re(w)}(a+cm)^{-\re(s)}\\
\leq & e^{\pi(|\im(s)|+|\im(w)|)} \sum_{c=1}^\infty\sum_{m=1}^\infty c\cdot c^{-k+\re(s)+\re(w)}(cm)^{-\re(s)}\\
=&e^{\pi(|\im(s)|+|\im(w)|)}\zeta(k-1-\re(w))\zeta(\re(s)),
\end{align*}
showing that the terms with $c>1$ sum to a holomorphic function on $\mathcal{F}$.

\noindent({\bf  2}) Next consider the terms with $c=1$. Separating the first term in the Hurwitz zeta functions, we have the following four terms
\begin{align*}
&\sum_{a=1}^\infty a^{-s}\int_0^1u^{s-1}(1-u)^{k-s-1}\cos(\pi(s+w)/2)(1-u/a)^{-w}du\\
+&\sum_{a=1}^\infty a^{-s}\int_0^1u^{s-1}(1-u)^{k-s-1}\cos(\pi(s-w)/2)(u/a)^{-w}du\\
+&\sum_{a=1}^\infty a^{-s}\int_0^1u^{s-1}(1-u)^{k-s-1}\cos(\pi(s+w)/2)\zeta\left(w,2-\frac{u}{a}\right)du
\\
+&\sum_{a=1}^\infty a^{-s}\int_0^1u^{s-1}(1-u)^{k-s-1}\cos(\pi(s-w)/2)\zeta\left(w,1+\frac{u}{a}\right)du.
\end{align*} The third and the fourth term are absolutely bounded by 
\[e^{\pi(|\im(s)|+|\im(w)|)}\zeta(\re(s))\zeta(\re(w)),\]
hence giving holomorphic functions on $\mathcal{F}$.
  Recognized as a beta integral, the second term is equal to   
\begin{align*}
&\frac{\Gamma(s-w)\Gamma(k-s)}{\Gamma(k-w)}\cos(\pi(s-w)/2)\zeta(s-w),
\end{align*}
which is meromorphic everywhere. Finally, we employ the elementary inequality for the first term
\[|(1-u/a)^{-w}|\leq (1-u)^{-\re(w)/a}, \quad a\in\mathbb{Z}_{>0},\re(w)>0, u\in (0,1).\]
The $a=1$ term gives
\[\frac{\Gamma(s)\Gamma(k-s-w)}{\Gamma(k-w)}\cos(\pi(s+w)/2),\]
where the rest in the first term gives a series absolutely convergent on $\mathcal{F}$, where it is bounded by
\[e^{\pi(|\im(s)|+|\im(w)|)}\zeta(\re(s)).\]

\noindent  ({\bf 3})    Putting everything together, on $\mathcal{F}$ we have, 
\begin{align*}
&A(s,w)\\
=&2(-1)^{\frac{k}{2}}\frac{(2\pi)^{k-w}\Gamma(w)}{\Gamma(s)\Gamma(k-s)\Gamma(k-w)}\\
&\times \left(\Gamma(s-w)\Gamma(k-s)\cos(\pi(s-w)/2)\zeta(s-w)+\Gamma(s)\Gamma(k-s-w)\cos(\pi(s+w)/2)\right)\\
&+ R(s,w),
\end{align*}
where $R(s,w)$ is holomorphic on $\mathcal{F}$ and is bounded by
\begin{equation}\label{Remainder}
2e^{\pi(|\im(s)|+|\im(w)|)}\frac{(2\pi)^{k-\re(w)}|\Gamma(w)|\zeta(\re(s))}{|\Gamma(s)\Gamma(k-s)|}\left(\zeta(k-1-\re(w))+\zeta(\re(w))+1\right).
\end{equation}
This concludes the result.
{\qed}

For fixed positive real numbers $\delta, T$, we consider the following smaller region   
$$\mathcal{F}_1:=\left\{(s,w)\in \mathbb{C}^2 : 
\frac{k-1}{2}<\re(s),\re(w)<\frac{k+1}{2}, \re(s)+\re(w)<k-\delta,  \im(s),\im(w)\in [-T,  T]\right\}. $$ Let us first make the expression in Proposition \ref{F} more symmetric.

\begin{lem}\label{F1}
We have the following identity on $\mathcal{F}_1$
\begin{align*} 
&\frac{2^{2}\pi^{k+1}\Gamma(k-1)}{\Gamma(s)\Gamma(w)\Gamma(k-s)\Gamma(k-w)}\sum_{f\in \mathcal{H}_k}\frac{L^*(f,s)L^*(f,w)}{\langle f,f\rangle}\\
=&\frac{(2\pi)^{s+w}}{\Gamma(s)\Gamma(w)}\zeta(k-s-w+1)+(-1)^\frac{k}{2}\frac{(2\pi)^{k-s+w}}{\Gamma(w)\Gamma(k-s)}\zeta(s-w+1)\\
&+(-1)^{\frac{k}{2}}\frac{(2\pi)^{k+s-w}}{\Gamma(s)\Gamma(k-w)}\zeta(w-s+1)+\frac{(2\pi)^{2k-s-w}}{\Gamma(k-s)\Gamma(k-w)}\zeta(s+w-k+1)\\
&+ R(s,w),
\end{align*}
where $R(s,w)$ (different from that in Proposition \ref{F}) is holomorphic on $\mathcal{F}_1$ and bounded by
\[\frac{(2\pi)^{k}}{|\Gamma(s)\Gamma(k-s)|}\]
up to a constant depending only on $T$ and $\delta$.
\end{lem}
{\pf}of Lemma \ref{F1}:
It follows from Proposition \ref{F} by simplifying the corresponding terms of the right-hand side therein. Explicitly, multiply both sides by $(2\pi)^{w}/\Gamma(w)$ and apply the functional equation of Riemann zeta function
\[\zeta(1-s)=2(2\pi)^{-s}\cos(\pi s/2)\Gamma(s)\zeta(s)\]
and then we have the desired third term. Recall that for $a,b\in\mathbb{C}$ and $t>0$, as $z\rightarrow\infty$ in the sector $|\textrm{arg}(z)|\leq \pi-t$,
\begin{equation}\label{Gamma}
\frac{\Gamma(z+a)}{\Gamma(z+b)}\sim z^{a-b}
\end{equation}
(see (5.11.12) of \cite{O}).
It follows that 
\[\frac{\Gamma(s)\Gamma(k-s)}{\Gamma(k-s)\Gamma(k-w)}=\frac{\Gamma(s)}{\Gamma(k-w)}=\frac{\Gamma(\frac{k}{2}-(\frac{k}{2}-s))}{\Gamma(\frac{k}{2}+(\frac{k}{2}+w))}\sim (k/2)^{s+w-k},\]
as $k\rightarrow\infty$, so the fourth term in Proposition \ref{F} and that of the present lemma can both be put into the remaining term $R(s,w)$.  Therefore, we can  make the expression symmetric by adding the zeta value in the fourth term. 
{\qed}

\medskip

\subsection{Proof of Theorem \ref{main-thm}}   Now we can prove the main theorem. Fix positive real numbers $T$ and $\delta$ and consider 
\[\mathcal{R}_{T,\delta}'=\{(s,w)\in\mathcal{F}_1\colon \re(s),\re(w)\leq k/2-\delta\}.\]
We prove instead that on $\mathcal{R}_{T,\delta}'$ the right-hand side of Lemma \ref{F1} is non-vanishing. To this end, we show that the first term is the main term, that is, we prove that when multiplied by $\Gamma(s)\Gamma(w)(2\pi)^{-s-w}$ the sum of the remaining terms  has limit $0$. 

As indicated in the proof of Lemma \ref{F1}, the fourth term therein can be put into the remaining term $R(s,w)$, while 
\[\left|\Gamma(s)\Gamma(w)(2\pi)^{-s-w}R(s,w)\right|\ll_{T,\delta}\left|\frac{(2\pi)^{k-s-w}\Gamma(w)}{\Gamma(k-s)}\right|.\]
By \eqref{Gamma},
\[\frac{(2\pi)^{k-s-w}\Gamma(w)}{\Gamma(k-s)}\sim (4\pi/k)^{k-s-w},\quad k\rightarrow\infty.\]
Since $\re(s)+\re(w)<k-\delta$, it follows that the fourth term and the remaining term give limit $0$ as expected.

So we only have to show that 
\begin{align*}
&\frac{(2\pi)^{k-2s}\Gamma(s)}{\Gamma(k-s)}\zeta(s-w+1)+\frac{(2\pi)^{k-2w}\Gamma(w)}{\Gamma(k-w)}\zeta(w-s+1)
\end{align*}
approaches $0$ as $k\rightarrow\infty$. Since the sum is holomorphic, the singularities of the two terms on $s=w$ cancel. Rewriting it as
\begin{align}\label{sing}
&\frac{(2\pi)^{k-2s}\Gamma(s)}{\Gamma(k-s)}\left(\zeta(s-w+1)-\frac{1}{s-w}\right)+\frac{(2\pi)^{k-2w}\Gamma(w)}{\Gamma(k-w)}\left(\zeta(w-s+1)-\frac{1}{w-s}\right)\\
&+\frac{(2\pi)^{k-2s}\Gamma(s)}{\Gamma(k-s)}\frac{1}{s-w}+\frac{(2\pi)^{k-2w}\Gamma(w)}{\Gamma(k-w)}\frac{1}{w-s}.\nonumber
\end{align}
The first two terms of \eqref{sing} approach $0$ by \eqref{Gamma} and the same argument as above, so it suffices to show
\[\frac{(2\pi)^{k-2s}\Gamma(s)}{\Gamma(k-s)}\frac{1}{s-w}+\frac{(2\pi)^{k-2w}\Gamma(w)}{\Gamma(k-w)}\frac{1}{w-s}\rightarrow 0,\quad k\rightarrow\infty.\]
It is clear that such a function is absolutely bounded by $\sup_{z}|g'(z)|$ where 
\[g(z)=\frac{(2\pi)^{k-2z}\Gamma(z)}{\Gamma(k-z)}\]
and $z$ ranges over $\{z\colon -1/2\leq \re(z)-k/2\leq -\delta, |\im(z)|\leq T\}$.
It is easy to obtain its derivative
\[g'(z)=\frac{(2\pi)^{k-2z}\Gamma(z)}{\Gamma(k-z)}\left(-2\log(2\pi)+\psi(z)-\psi(k-z)\right),\]
where $\psi(z)=\Gamma'(z)/\Gamma(z)$.
Since $\psi(z)\sim \log z$ (see (5.11.2) of \cite{O}), we have $\sup_{z}|g'(z)|\rightarrow 0$ as desired.

The above argument shows that on $\mathcal{R}_{T,\delta}'$, if $k$ is large enough,
\[\sum_{f\in \mathcal{H}_k}\frac{L^*(f,s)L^*(f,w)}{\langle f,f\rangle}\neq 0.\]
Finally, by applying the functional equations of $E_{s,k-s}^*(z,w)$ we obtain the theorem. {\qed}

\medskip
\subsection{Proof of Theorem \ref{without}}
By Theorem \ref{main-thm} and the symmetries, we only have to show that for each $(s_0,w_0)=(1/2+it_0,w_0)$ with $t_0\neq 0, \re(w_0)\leq 1/2-\delta$ and $\im(w_0)\neq 0$, there exists a neighborhood on which the right-hand side of the equation of Lemma \ref{F1} (after shifting by $\frac{k-1}{2}$) is non-vanishing when $k$ is large. Following the proof of Theorem \ref{main-thm}, we see that the remainder on the right-hand side of the equation of Lemma \ref{F1} is dominated by the first two terms as $k\rightarrow\infty$. Therefore, we need to prove that on some neighborhood of $(s_0,w_0)$, the quantity
\begin{align*}
&1+(-1)^{\frac{k}{2}}(2\pi)^{1-2s}\frac{\Gamma(s+\frac{k-1}{2})}{\Gamma(\frac{k+1}{2}-s)}\frac{\zeta(1+s-w)}{\zeta(2-s-w)}+(-1)^{\frac{k}{2}}(2\pi)^{1-2w}\frac{\Gamma(w+\frac{k-1}{2})}{\Gamma(\frac{k+1}{2}-w)}\frac{\zeta(1+w-w)}{\zeta(2-s-w)}\\&+(2\pi)^{2-2s-2w}\frac{\Gamma(s+\frac{k-1}{2})}{\Gamma(\frac{k+1}{2}-s)}\frac{\Gamma(w+\frac{k-1}{2})}{\Gamma(\frac{k+1}{2}-w)}\frac{\zeta(s+w)}{\zeta(2-s-w)}
\end{align*}
is non-vanishing when $k$ is large. 

The third and the fourth term approach $0$ when $k$ approaches $\infty$, as we have seen in the proof of Theorem \ref{main-thm}. Now \[\frac{\Gamma(s+\frac{k-1}{2})}{\Gamma(\frac{k+1}{2}-s)}=\left(\frac{k}{2}\right)^{1-2s}\left(1+O(k^{-1})\right), k\rightarrow\infty,\]
so together with the Assumption, it implies the existence of a desired neighborhood, on which the sum of first two terms stays away from $0$.
{\qed}

%%%%%%%%%%%%%%%%%%%%%%%%%%%%%
\bibliographystyle{amsplain}

\end{document}